\documentclass[journal, 11 pt, onecolumn]{IEEEtran}
\makeatletter

\let\proof\@undefined
\let\endproof\@undefined
\makeatother

\usepackage{amsthm,amssymb,euscript,amscd,cite,array,color}
\usepackage[cmex10]{amsmath}
\usepackage{epsfig}

\usepackage{graphicx}

\usepackage[tight,footnotesize]{subfigure}
\newtheorem{lemma}{Lemma}

\newtheorem{prop}{Proposition}
\newtheorem{remark}{Remark \textrm}

\newcommand{\pr}{\textnormal{Pr}}
\newcommand{\E}{\mathbf{E}}

\newcommand{\disp}{\displaystyle}

\newcommand{\sign}{ \mbox{sign}}

\newcommand{\calF}{{\mathcal{F}}}

\newcommand{\ba}{\begin{array}}
\newcommand{\ea}{\end{array}}

\newcommand{\ee}{\end{equation}}
\newcommand{\be}{\begin{equation}}

\newcommand{\vsp}{\vspace{0pt}}

\newcommand{\Erdos}{Erd\"{o}s}
\newcommand{\Renyi}{R\'{e}nyi}
\newcommand{\Czirok}{Czir\'{o}k}

\newcommand{\tabfrac}[2]{%
	\setlength{\fboxrule}{0pt}%
	\fbox{$ \displaystyle \frac{#1}{#2}$}%
}

\begin{document}

\title{Phase Transitions on Fixed Connected Graphs and Random Graphs in
the Presence of Noise}


\author{Jialing Liu, Vikas Yadav, Hullas Sehgal, Joshua M. Olson, Haifeng Liu, and Nicola Elia 
\thanks{J. Liu was partially supported by NSF under Grant
ECS-0093950.    A preliminary version of this paper has appeared
in Proceedings of the 44th IEEE Conference on Decision and Control and
European Control Conference (CDC-ECC'05).  This work was performed when the authors were with the Department of Electrical and Computer Engineering, Iowa State
University, Ames, IA 50011, USA.}
\thanks{ J. Liu is with Motorola Inc., 600 N. US-45, Libertyville, IL 60048 USA (e-mail: jialingliu@motorola.com).}
\thanks{ V. Yadav is with Garmin International, Olathe, KS 66062 USA (e-mail: vikas.yadav@garmin.com).}
\thanks{ H. Sehgal is with Electrical Engineering Department,
University of Minnesota (Twin Cities), Minneapolis, MN 55455 USA (e-mail: sehga008 @umn.edu).}
\thanks{ J. M. Olson is with Raytheon Missile Systems, Tucson, AZ 85743 USA (email: joshua\_m\_olson@raytheon.com). }
\thanks{H. Liu is with California Independent System Operator, Folsom, CA 95630 USA (e-mail: hliu@caiso.com).}
\thanks{N. Elia is with the Department of
Electrical and Computer Engineering, Iowa State University, Ames,
IA 50011 USA (e-mail: nelia@iastate.edu).}
}

 \maketitle 

\begin{abstract} In this paper, we study the phase transition behavior emerging from
the interactions among multiple agents in the presence of noise.
We propose a simple discrete-time model in which a group of
non-mobile agents form either a fixed connected graph or a random
graph process, and each agent, taking bipolar value either $+1$ or
$-1$, updates its value according to its previous value and the
noisy measurements of the values of the agents connected to it.  We present
proofs for the occurrence of the following phase transition
behavior: At a noise level higher than some threshold, the system
generates symmetric behavior (vapor or melt of magnetization) or disagreement; whereas at a noise level lower than
the threshold, the system exhibits spontaneous symmetry breaking (solid or magnetization) or consensus.
The threshold is found analytically.  The phase transition 
occurs for any dimension.  Finally, we demonstrate the phase transition behavior
and all analytic results using
simulations. This result may be found useful in the study of the
collective behavior of complex systems under communication
constraints.
\end{abstract}


\section{Introduction}

A phase transition in a system refers to the sudden change of a 
system property as some
parameter of the system crosses a certain 
threshold value.  Phase transitions have
been observed in a wide 
variety of studies, such as in physics, chemistry, biology,
complex systems, computer science, and random graphs, to list a 
few.  It leads to long
term attention in the literature, from 
physicists such as Ising \cite{Ising87} in the
1920's to 
mathematicians such as \Erdos~ and \Renyi~\cite{erdos60} in the 1960's, from
complex systems theorists such as Langton \cite{langton90} in the 1990's to control
scientists such as Olfati-Saber \cite{Olf05} in the 2000's.

Ising and other physicists have thoroughly studied the simple but
``realistic enough" Ising model, for the understanding of phase
transitions in magnetism, lattice gases, etc.  In an Ising model,
each node can take one of two values, and the neighboring nodes have an
energetic preference to take the same value, under some
constraints such as a temperature one.  It is observed that, for
an Ising model with dimension at least 2, a temperature higher
than a critical point leads to symmetric behavior (e.g., ``melt"
of magnetization, or vapor), whereas a temperature lower than that
point leads to asymmetric behavior (e.g., magnetization, or
liquid).  The Ising model is a discrete-time discrete-state model,
and is closely related to the Hopfield networks and cellular
automata.

\Erdos~and \Renyi~\cite{erdos60} showed that, graphs of sizes
slightly less than a certain threshold are very unlikely to have
some properties, whereas graphs with a few more edges are almost
certain to have these properties. This is called a phase transition
of random graphs, see for example \cite{janson00}.

Viscek \emph{et al} \cite{Vic&Czi95} showed that a two-dimensional
nonlinear model exhibits a phase transition in the sense of
spontaneous symmetry breaking as the noise level crosses a
threshold. This model consists of a two-dimensional square-shaped
box filled with particles represented as point objects in
continuous motion.  The following assumptions are also adopted: 1)
the particles are randomly distributed over the box initially; 2)
all particles have the same absolute value of velocity; and 3) the
initial headings of the particles are randomly distributed.  Each
particle updates its heading using the average of its own heading
and the headings of all other particles within a radius $r$, which
is called the \emph{nearest neighbor rule}. Included for each particle in this model
is a random noise (which may be viewed as measurement noise or actuation noise) with a uniform probability distribution on the
interval $[-\eta, \eta]$. The result of \cite{Vic&Czi95} is to
demonstrate using simulations that a phase transition occurs when
the noise level crosses a threshold which depends on the particle
density.  Below the threshold, all particles tend to align their
headings along some direction, and above the threshold, the
particles move towards different directions.  \Czirok~\emph{et al}
\cite{Czi&Vic99} presented a one-dimensional model which also
exhibits a phase transition for a group of mobile particles.  These
two models are discrete-time continuous-state models.

Schweitzer~\emph{et al} \cite{schweitzer} studied the
spatial-temporal evolution of
the populations of two species, where the update scheme depends nonlinearly on the
local frequency of species.   Depending on the probability of transition from one
species to the other, the system evolves to either extinction of one species
(agreement) or non-stationary co-existence or random co-existence (disagreement).

We note that phase transition problems are sometimes associated with flocking / swarming / formation / consensus / agreement problems.  Though the interest and focus of these problems are often independent of the phase transition study, these problems typically exhibit phase transitions when parameters, conditions, or structures change.  These problems and the phase transition problems may also share some common techniques in order to establish stability / instability over similar underlying models, such as common Lyapunov function techniques, graph theoretic techniques, and stochastic dynamical systems techniques.
More specifically, the phase transitions occurring in flocking may be classified into two categories: angular phase transitions that leads to alignment (see e.g. \cite{Vic&Czi95}), and spatial self-organization in which multiple agents tend to form special patterns or structures in space, such as lattice type structures. Examples of the latter category include \cite{Mog&Keshet1996, Lev&Rapp&Coh01, Olf06}.  In \cite{Mog&Keshet1996}, Mogilner and Edelstein-Keshet investigated swarming in which the dynamical objects interacts depending on angular orientations and spatial positions, and a phase transition is observed.  In \cite{Lev&Rapp&Coh01}, Levine \emph{et al} presented a simple model to
study spatial self-organization in flocking showing that all the
agents tend to localize in a special pattern in one- and in two-dimensions with all-to-all communication. 
We refer to \cite{OlfSab&Murray03, OlfSab&Murray04,Olf05, Olf06, Hatano&Das&Mesbahi05} for some recent studies of phase transitions and
the consensus / agreement problems over networks.  Olfati-Saber \cite{Olf05}
studied the consensus problem using a random rewiring algorithm (see also \cite{watts98}) to
connect nodes, and showed that the Laplacian spectrum of this network may undergo a
dramatic shift, which is referred to as a spectral phase transition and leads to
extremely fast convergence to the consensus value. In \cite{Hatano&MesbahiTAC05}, Hatano and Mesbahi established agreement of multiple agents over a network that forms an \Erdos~random graph process, in which each agent updates its state linearly according to the perfect state information of its nearest neighbors.  Hatano and Mesbahi also studied another facet of the distributed agreement problem in the stochastic setting in \cite{Hatano&Das&Mesbahi05}, namely the agreement over noisy network that forms a Poisson random graph.

Jadbabaie \emph{et al} \cite{Jad&Lin02} provided a rigorous proof
for the alignment of moving particles under the nearest neighbor rule without measurement noise or actuation noise.  Different from the switched nonlinear model used in
\cite{Vic&Czi95}, the model in
\cite{Jad&Lin02} is a switched linear model. Furthermore, this
model also assumes that over every finite period of time the
particles are \emph{jointly connected} for the length of the
entire interval. Due to the
noiseless assumption made in \cite{Jad&Lin02}, the phase
transitions observed in \cite{Vic&Czi95} will not occur. Under these assumptions, Jadbabaie \emph{et al}
proved that the nearest neighbor rule leads to alignment of all
particles. One may be interested in finding Lyapunov functions
(preferably quadratic) to show the convergence or alignment (see \cite{OlfSab&Murray03, OlfSab&Murray04} for convergence proofs based on common Lyapunov functions for models different from \cite{Jad&Lin02}).
However, \cite{Vic&Czi95} showed that a common quadratic Lyapunov
function does not exist for this switched linear model.  On the
other hand, a non-quadratic Lyapunov function can be constructed
to prove the convergence, as suggested by
Megretski~\cite{meg:lyap} and later independently found by Moreau
\cite{lucMoreau05}.  See also \cite{Jad&Tan1,Jad&Tan2} for
extension of \cite{Jad&Lin02}.

In this paper, we propose a discrete-time discrete-state model in
which a group of agents form either a fixed connected graph or a
random graph process, and each agent (node) updates its value
according to its previous value and the noisy measurements of the
neighboring agent values.  We prove that, when the noise level
crosses some threshold from above, the system exhibits spontaneous
symmetry breaking. We may view that the high noise level
corresponds to high temperature (or strong thermal agitation),
where the molecules exhibit disorder and symmetry; and the low
noise level corresponds to low temperature, where the molecules
exhibit order and asymmetry.

We emphasize that the proposed model is rather simple and hence
admits a complete mathematical analysis of the phase transition
behavior. First, the phase transition in a fixed connected graph
presented in this paper is simpler than the phase transition in
the Ising model.  As one indicator of the simplicity, note that the Ising model of dimension higher than two involves intractable computation complexity when attempting to solve for the value for each node under the temperature constraint, namely, such a problem is NP-complete \cite{cipra_npscience}.  Also
the Ising
model needs dimension two or higher to generate the phase transition,
whereas our phase transition occurs for any dimension. To the best of our knowledge,
the proposed model is one of the simplest that exhibits a phase transition
in a fixed graph, and is mathematically provable to generate a sharp phase
transition.   Note that many other phase transitions elude
rigorous mathematical analysis due to their complexity
\cite{langton90,Vic&Czi95,Czi&Vic99,cipra_npscience,Olf05}.
Second, the phase transition on a random graph is also simpler
than the phase transition on a random graph observed in
\cite{erdos60}. Compared with the models in \cite{Vic&Czi95} and
\cite{Czi&Vic99}, our models have discrete-states and do not allow
the mobility of agents, which greatly simplifies the systems
dynamics and allows rigorous proofs of the phase transition
behavior. The simplicity of our phase transitions may help us to
identify the essence of general phase transition phenomena.

Our study also sheds light on the research on the consensus
problems, cooperation of multiple-agent systems, and collective
behavior of complex systems, all under communication constraints.
Hence, this study fits into the general framework
of investigating the interactions between control/dynamical
systems and information; see e.g.
\cite{wolf:phd,mitter:talk,tati:phd,sahai:phd,elia_c5,fax04,liu:phd} and references
therein.  More specifically, we may interpret our phase transition in
the consensus problem framework, where the disagreement due to
unreliable communication is replaced by agreement when the
communication quality improves to a certain level.  In other words, our work characterizes the significance of information in reaching agreement.  However,
unlike the average-consensus problem (cf. \cite{OlfSab&Murray04})
with the properties that, 1) there exists an invariant quantity
during the evolution, and 2) the limiting behavior reaches the
average of the initial states of the system, our models reach agreement
without these properties when the noise level is low.  This is
because the presence of noise prevents the conservation of the sum
of the node values during the evolution.  The study of entropy flows (or
information flow) \cite{mitter:kalman,wolf:phd} may help identify
an invariant quantity of the system. We remark that a more thorough study of the
consensus problem raised in this paper is beyond the scope of this
paper and will be pursued elsewhere.

\textit{Organization:} In Section 2 we introduce the models. In Section 3 we state our
main results and provide the proofs. In Section 4 we present numerical examples.
Finally we conclude the paper and discuss future research directions.

\section{Models on the graphs} 

This section introduces some of the terms that are frequently used
in this paper as well as the two models to be investigated.  We
focus only on undirected graphs.

\subsection{Graphs and random graph processes}

A $\textit{graph}$ $G := (V, E)$ consists of a set $V:=
\{1,2,...,N\}$ of elements called vertices or nodes, and a set $E$ 
of node pairs called edges, with $E \subseteq E_c := \{(i,j)| i,j 
\in V\}$. Such a graph is $\textit{simple}$ if it has no self
loops, i.e. $(i,j) \not \in E$ if $ i=j$.  We consider simple 
graphs only.  A graph $G$ is $\textit{connected}$ if it has a path 
between each pair of distinct nodes $i$ and $j$, where by a
\emph{path} between nodes $i$ and $j$ we mean a sequence of
distinct edges of $G$ of the form $(i, k_1), (k_1, k_2), \ldots
,(k_m,j) \in E$. Radius $r$ from node $i$ to node $j$ means that 
the minimum path length, i.e., the minimum number of edges
connecting $i$ to $j$, is equal to $r$.

A $\textit{fixed}$ graph $G$ has a node set $V$ and an edge set
$E$ that consists of fixed edges, that is, the elements of
$E$ are deterministic and do not change dynamically with time.

A $\textit{random}$ graph $G$ consists of a node set $V$ and an
edge set $E:=E(\omega)$, where $\omega \in \Omega$,
$(\Omega,\calF,P)$ forms a probability space.  Here $\Omega$ is
the set of all possible graphs of total number of $n$, where
\be n:=2^{ \frac{ N(N-1) } { 2 } }; \label{totngraph} \ee
$\calF$ is the power set of $\Omega$; and $P$
is a probability measure that assigns a probability to every
$\omega \in \Omega$. In this paper, we focus on the well-known
\Erdos~random graphs \cite{janson00}, namely, it holds that
\be P(\omega)= \frac{1}{n}. \label{Pomega}\ee
In other words, we can view each $E(\omega)$ as a result of
$N(N-1)/2$ independent tosses of fair coins, where a head
corresponds to switching on the associated edge and a tail
corresponds to switching off the associated edge. Notice that the
introduction of randomness to a graph implies that, all results for
random graphs hold asymptotically and in a probability
sense, such as ``hold with probability one".

A \emph{random graph process} is a stochastic process that
describes a random graph evolving with time. In other words, it is
a sequence $\{G(k)\}_{k=0}^\infty$ of random graphs (defined on a
common probability space $(\Omega,\calF,P)$) where $k$ is
interpreted as the time index (cf. \cite{janson00}).  For a random
graph process, the edge set changes with $k$, and we denote the
edge set at time $k$ as $E(k)$. In this paper, we assume that the
edge formation at time $k$ is independent of that at time $l$, if
$k \neq l$.

The \emph{neighborhood} $N_i(k)$ of the $i$th node at time $k$ is
a set consisting of all nodes within radius 1, including the $i$th
node itself.  The value that a node assumes is its \emph{node
value}. The \emph{valence} or \emph{degree} of the $i$th node is
$(|N_i(k)|-1)$, where $|N_i(k)|$ denotes the number of elements in
$N_i(k)$. The \emph{adjacency matrix} of $G(k)$ is an $N \times N$
matrix whose $(i,j)$th entry is $1$ if the node pair $(i, j) \in
E(k)$ and 0 otherwise.  Note that the graphs can model lattice
systems with any dimension.

\subsection{System on a graph}

A \emph{system on a graph} consists of a graph, fixed or forming a
random process, an initial condition that assigns each node (agent) a node
value, and an update rule of the node values.  In this paper, we
assume that each node can take value either $+1$ or $-1$, and the
$\textit{update rule}$ for the $i$th node at the $(k+1)$st instant
is given by
\be x_i(k+1) = \textnormal{sign} \left[v_i(k) + \xi_i(k) \right],
\label{eq:update}\ee
where $\xi_i(k)$ is the \emph{noise} random variable, uniformly
distributed in interval $[-\eta,\eta]$ and independent across time
and space and independent of the initial condition $x(0)$, and
\be v_i(k) := \frac{\sum_{j \in N_i(k)} x_j(k)}{|N_i(k)|} ;\ee
that is, $v_i(k)$ is the average of the node values in the
neighborhood $N_i(k)$.  Here $\eta$ is called the \emph{noise
level}.  This update rule resembles the one in \cite{Czi&Vic99},
with their antisymmetric function being replaced by a sign
function.  It may also be viewed as a specific update rule for a
Hopfield neuron whose connections with others are noisy.

The \emph{state of the system} at time instant $k$, denoted
$x(k)$, is the collection of all node values $(x_1(k),\cdots,x_N(k))$.  The \emph{state sum} at time instant $k$, denoted
$S(k)$, is defined as
\be S(k) := \sum_{i=1} ^N x_i (k) . \ee
With a slight abuse of notation, we represent the state with all +1s and all -1s (i.e. the consensus states) as $+N$ and $-N$, respectively.
We call a state \emph{transient} if this state reappears
with probability strictly less than one. We call a state
\emph{recurrent} if this state reappears with probability one. We
call a state $X$ \emph{absorbing} if the one-step transition
probability from $X$ to $X$ is one.

\subsection{Model with a fixed graph}

The first model considered is a system on a fixed graph.  In this model, the node
connections or the edges remain unchanged throughout. Hence, every node has a fixed
neighborhood at all times, and the degree of each node as well as the adjacency matrix
are constant. The node value gets updated according to the update rule
(\ref{eq:update}). We will assume that the fixed graph is connected.  An example of
such a fixed graph model is a communication network with fixed nodes and fixed but
noisy channels.  Another example is a Hopfield network with fixed neurons and fixed
but noisy channels connecting them.

\subsection{Model with a random graph process}
The second model considered is a system on a graph forming a 
random process.  In this model, the node connections, namely the 
edges of the random graph, change dynamically throughout, and the 
edge formations at time $k$ are random according to distribution
$P(k)$. Hence every node may have different neighborhoods at 
different times, and the adjacency matrix and degrees change with 
time. The node value gets updated also according to the update 
rule (\ref{eq:update}).  An example of this model could be an
ad-hoc sensor network in which the communication links between the
sensors appear and disappear dynamically.  Another example is an erasure network in 
which the communication channels are noisy and erasing with some
probability, see for example \cite{julian_erasure02}.

In both models, the system state can take $2J$ values, where
\be J := 2^{ N - 1 } \label{ss_size} \ee
and the state sum takes values in the set
$\mathcal{N}:=\{-N,-N+2,\cdots,N-2,N\}$, where $N \geq 2$ is the total number of nodes.  Note that $|\mathcal{N}|=N+1 \geq 3$. Both models also form Markov chains,
since the next state does not depend on previous state if the current state is given.

We use $\xi(k)$ to represent $(\xi_1(k),\cdots,\xi_N(k))$, $\xi^k$ to represent
$(\xi(0),\cdots,\xi(k))$, $G^k$ to represent $(G(0),\cdots,G(k))$, and $x(k)$ to represent $(x_1(k),\cdots,x_N(k))$.  \vsp

\section{Main results and proofs}

Our main result states that,\emph{ for a system on a fixed 
connected graph or on a graph forming a random process, there is a
provable sharp phase transition when the noise level crosses some
threshold}. Here the phase transition is in the sense that the
symmetry exhibited at high noise level is broken suddenly when the
noise level crosses the threshold from above, or equivalently the
disagreement (or disorder) of the nodes at high noise level becomes agreement (or order)
below the threshold. In what follows, we first discuss the case in
which the graph has a fixed structure, and then the case in which
the graph forms a random process.

\subsection{Model with a fixed graph}

\begin{prop} \label{prop:fix}
For any given fixed connected graph, let $D$ be the maximum number
of nodes in one neighborhood.

i) If the noise level is such that $\eta \in (1-2/D, 1]$, then the
system will converge to \emph{agreement}, namely all nodes will
converge to either all $+1$s or all $-1$s.

ii) If the noise level is such that $\eta > 1$, then $\E
 S(k)$ tends to zero as $k$ goes to infinity,
i.e., the system will converge to \emph{disagreement} in which
approximately half of the nodes are $+1$s and the other half are
$-1$s.

\end{prop}

\begin{remark} \rm Notice that $(1-2/D)$ is guaranteed to be
nonnegative for any connected graph with more than one node, since
$D\geq 2$. Note also that if $\eta <(1-2/D)$, the system does not
necessarily converge to states $\pm N$. To see this, simply
consider a one-dimensional cellular automaton with $N$ nodes
forming a circle.  The neighborhood of a node is defined as one
node to the left, one node to the right, and itself. Therefore
$D=3$, and if $\eta< 1/3$, the update rule becomes a majority
voting rule. Then the initial condition $x(0)$ of the system with
alternate $+1$s and $-1$s will lead to constant oscillations
between $x(0)$ and a left cyclic shift of $x(0)$, i.e., it will not reach
agreement if $\eta < 1/3$. However, this does not mean that in
general our condition $1 \geq \eta > (1 - 2/D)$ is a necessary
condition for agreement; a sufficient and necessary condition is
under current investigation. Attractors like this $x(0)$ may be
viewed as local attractors (whereas $\pm N$ may be viewed as
global attractors) which can be eliminated by considering a
\emph{randomized} graph, see the next subsection.
\end{remark}

The proof of Proposition \ref{prop:fix} needs the following lemmas.

\begin{lemma} \label{lemma:fix1}
For any given fixed connected graph, if $\eta \in (1-2/D, 1]$,
then the states $\pm N$ are absorbing, and all other states are
transient.
\end{lemma}

\begin{lemma}\label{lemma:fix2}
For any given fixed connected graph, if $\eta > 1$, then the
states form an ergodic Markov chain with a unique steady-state
distribution for any initial condition $x(0)$.
\end{lemma}

\textbf{Proof of Lemma \ref{lemma:fix1}:} At states
$\pm N$, the noise is not strong enough to flip any node value. 
Thus, $\pm N$ are absorbing.  On the other hand, all other states are
neither absorbing nor recurrent.  To see this, let $M \neq \pm N$
be any state, which leads to that $M$ contains a mixture of $+1$s and
$-1$s.  Due to the connectivity of the graph, we can always find a node $i$ with node value
$x_i(k)=-1$ whose neighborhood $N_i(k)$ (including $x_i(k)$
itself) contains both $+1$s and $-1$s.  Then for such $x_i(k)$, it holds that
\be \left| v_i(k) \right| \leq \frac{D-2}{D} ,\ee
with equality if only one node in $N_i(k)$ has a different sign
than all other nodes and if $N_i(k)$ contains $D$ nodes. Hence a
noise larger than $(D-2)/D$ flips $x_i(k)$.  Precisely,
\be \ba{lll} &&\Pr [ x_i(k+1)=+1|x_i(k)=-1 ] \\
&=& \disp \Pr \left[
v_i(k) +\xi_i(k) >0 |x_i(k)=-1 \right] \\
& \geq &  \disp  \Pr \left[ \left. \xi_i(k)
>  \disp \frac{D-2}{D} \right| x_i(k)=-1 \right] \\
&=&  \disp \frac{1}{2} \left(1-\frac{D-2}{D\eta} \right)  > 0. \ea \ee
Note that the conditioning is removed due to the independence assumptions on noise.  Thus, for state $M$, the probability that only $x_i$ flips and no other node changes
its value is non-zero. This follows that, with a positive probability the state sum for $M$ will be increased by $2$. Likewise, with a positive probability $M$ can be decreased by 2.
Since $M \neq \pm N$ is an arbitrary state, by induction, the probability of
transition (in possibly multiple steps) from $M$ to $\pm N$ is nonzero. So $M$ is
transient. \endproof  \vsp \vsp

\textbf{Proof of Lemma \ref{lemma:fix2}:} It is sufficient to
prove that the state forms an irreducible and aperiodic Markov
chain.

To see the irreducibility, note that if $\eta
>1$,  $M \neq \pm N$ can jump to any other states with a positive probability, similar to Lemma
\ref{lemma:fix1}. Additionally, $\pm N$ can also jump to any other states with a positive probability.  For
state $+N$, it holds that
\be \ba{lll}
 &&\Pr [ x_i(k+1)=-1|x_l(k)=+1, l=1,\cdots,N ] \\
&=& \disp \Pr \left[ \left. v_i(k) + \xi_i(k)
<0 \right| x_l(k)=+1, l=1,\cdots,N   \right] \\
&=& \disp \Pr [ \xi_i(k) <-1 |x_l(k)=+1, l=1,\cdots,N ]\\
&=& \disp \frac{1}{2\eta} (\eta-1) >0, 
\ea \ee
so any node can flip its value with a positive probability.
Similar result holds for state $-N$. Then this Markov chain is irreducible.

To see the aperiodicity, let us use $-{x}_i$ to denote the flipped $x_i$.  The state transition cycle from $(x_1 (k), x_2(k), *)$ to $(-{x}_1 (k), -x_2(k), *)$ to $(-{x}_1 (k), x_2(k), *)$ and back to $(x_1 (k), x_2(k), *)$ has period 3, where $*$ is any fixed configuration for $(x_3(k), \cdots, x_N(k))$.  However, the state transition cycle from $(x_1 (k), \Delta)$ to $(-{x}_1 (k), \Delta)$ and back to $({x}_1 (k), \Delta)$ has period 2, where $\Delta$ is any fixed configuration for $(x_2(k), \cdots, x_N(k))$.
Note that such cycles occur with positive probabilities.  Then the Markov
chain is aperiodic. \endproof  \vsp \vsp

\textbf{Proof of Proposition \ref{prop:fix}:} If $\eta \in (1-2/D,
1]$, by Lemma \ref{lemma:fix1}, the associated Markov chain will
converge to either $+N$ or $-N$ with probability 1, namely
agreement.  If $\eta >1$, from Lemma \ref{lemma:fix2} we know that
the associated Markov chain is ergodic, and notice that the Markov
chain has a symmetric structure for states $x$ and $-x$.  Then $
\pi(x) = \pi(-x)$ (rigorous proof is included in Appendix), where $\pi(x)$ is the stationary probability
of state $x$. Hence the expectation of state sum is
\be \E _{ S \sim \pi } S = \sum_{x} \left( \pi(x) \sum _{i=1}^N x_i \right) =0.\ee
Therefore, $\E  S(k)$ converges to zero, and the numbers of $+1$s and
$-1$s will asymptotically become equal. \vsp \endproof

\vsp \vsp
\subsection{Model with a random graph process} \vsp

For an \Erdos~random graph, we assume that the edge connections are randomly and
independently changing from time to time.  The randomization of the connections
symmetrizes the system behavior and leads to agreement even for an arbitrarily small but
positive noise level.

\begin{prop} \label{prop:random}
Consider an \Erdos~random graph process.

i) If the noise level is such that $0<\eta \leq 1$, then the
system will converge to \emph{agreement}, namely the state will
converge to $+N$ or $-N$.

ii) If the noise level is such that $\eta > 1$, then $\E  S(k) $ exponentially converges to zero
with decay exponent $\ln \eta$ as $k$ goes to infinity, i.e., the
system will exponentially converge to \emph{disagreement} in which
about half of the node values are $+1$s and the other half are
$-1$s.

\end{prop}

The proof of
this proposition needs the following lemmas.  We remark that it is
straightforward to generalize the lemmas to a binomial random
graph, in which the probability of forming an edge is changed from
$0.5$ to an arbitrary $p \in (0,1)$.

\begin{lemma} \label{lemma:random1}
For any \Erdos~random graph process, if $\eta \in (0, 1]$, then $\pm
N$ are absorbing, and all other states are transient.
\end{lemma}

\begin{lemma} \label{lemma:random2}
For any \Erdos~random graph process, if $\eta > 1$, then it holds
that $\E S(k)$ exponentially tends to
zero as $k$ goes to infinity.  The decay exponent is $\ln \eta$.
\end{lemma}

\textbf{Proof of Lemma \ref{lemma:random1}:} If $0< \eta \leq 1$,
it is easy to see that $\pm N$ are absorbing. For any state $M
\neq \pm N$, it holds that $M$ must be a mixture of both $+1$s and
$-1$s. Hence we can find $i$ and $j$ in $V$ such that $x_i(k)=-1$ and $x_j(k)=+1$.
Since each of the $n$ graphs (recall (\ref{totngraph})) has a positive probability, the
probability that $x_i$ is connected to $x_j$ only is positive.
Then in this case, $v_i(k)$ is 0 and hence an
arbitrarily small but positive noise may flip $x_i$ with a
positive probability. In addition, each node other than $x_i$ has
a positive probability to keep its previous value, thus with a
positive probability, the state sum for $M$ can be increased by $2$.  Therefore any $M \neq
\pm N$ are transient. 

\endproof \vsp

\textbf{Proof of Lemma \ref{lemma:random2}:} For any \Erdos~random
graph, if $\eta>1$, then no state is absorbing, since with
a positive probability the noise can flip any node value in any
configuration.  Therefore, with a nonzero probability the state of
the system can jump to any other states.

Now let us analyze the evolution of $\E S(k)$.  Fix the time to be $k$.
Assume $x(k)$ is given.  Then for each $i$, $x_i(k+1)$ is given
by (\ref{eq:update}).  The randomness in $x_i(k+1)$ is due to the
noise $\xi_i(k)$ and the graph $G(k)$. It holds that
\be \ba{lll} &&\E  [ x_i(k+1) | x(k) ] \\
&=& \displaystyle \E \:
\sign [ v_i(k) + \xi_i (k)  | x(k) ] \\
&=& \displaystyle \pr [v_i(k) + \xi_i (k) >0 | x(k) ] \times (+1) \\ &&
\;+ \pr [v_i(k) + \xi_i
(k) <0 | x(k) ] \times (-1) \\
&=& \displaystyle \pr [ \xi_i (k) > - v_i(k)  | x(k) ] - \pr [ \xi_i (k) < - v_i(k)  | x(k) ] \\
&=& \displaystyle \sum_{v_i(k)} \pr [ \xi_i (k) > - v_i(k) | v_i(k)
] \pr[v_i(k) | x(k) ]
\\  && \displaystyle\;- \sum_{v_i(k)} \pr [ \xi_i (k) < - v_i(k) | v_i(k) ] \pr[v_i(k) | x(k) ] \\
&=& \displaystyle \sum_{v_i(k)} \left[ \frac{\eta + v_i(k)}{2\eta}
 \displaystyle -  \frac{\eta - v_i(k)}{2\eta} \right] \pr[v_i(k) | x(k) ] \\
&=& \displaystyle \sum_{v_i(k)}  \frac{v_i(k)}{\eta}  \pr[v_i(k) | x(k) ]
\\ & =& \displaystyle \frac{1}{\eta} \E [ v_i(k) | x(k) ]. \ea \label{Exik} \ee

Then we compute $\E ( v_i(k) | x(k) )$.  Since conditioned on $x(k)$, the randomness in $v_i(k)$ comes from $G(k)$ only, this expectation boils down to the expectation of the
average of node values in a neighborhood, averaged over all possible $n$
graph structures.  Let us count in the $n$ graph structures the
number of different neighborhood types containing node $i$. Among
those neighborhoods containing node $i$, there are
\be  \bar{m} := 2^{(N-1)(N-2)/2} \times {N-1 \choose m}\ee
types of neighborhoods for which $|N_i(k)| = (m+1)$ where
$m=0,1,\cdots,N-1$.  To see this, simply notice that the graph formed by nodes other than $i$ can have any edge formation and hence the number of types of $2^{(N-1)(N-2)/2}$, and that node $i$ needs to select $m$ out of the other $(N-1)$ nodes in order to have $|N_i(k)| = (m+1)$.

Therefore,
\be \ba{lll} &\E [ v_i(k) | x(k) ] \\
= & \displaystyle \sum_{G(k)} [v_i(k)| x(k), G(k)] \pr[G(k)]\\
= & \displaystyle \frac{1}{n} \sum_{G(k)} \left[ \left. \frac{\sum_{j\in N_i(k)} x_j(k)}{|N_i(k)|} \right| x(k), G(k) \right]  \\
= & \displaystyle \frac{1}{n} \sum_{G(k)} \left[ \left. \frac{x_i(k)}{|N_i(k)|} \right| x(k), G(k) \right]  + \\
 & \quad \displaystyle \frac{1}{n} \sum_{G(k)} \left[ \left. \frac{\sum_{j\in N_i(k),j \neq i} x_j(k)}{|N_i(k)|}  \right| x(k), G(k) \right] .\ea \ee
Now first note that
\be  \disp  \sum_{G(k)} \left[ \left. \frac{x_i(k)}{|N_i(k)|} \right| x(k), G(k) \right] 
=   \disp  \sum_{m = 0} ^ {N-1} \left[ \left. \frac{x_i(k)} {m+1} \bar{m} \right| x(k) \right]. \ee
Then note that in the summation
\be \sum_{G(k)} \left[ \left. \frac{\sum_{j\in N_i(k),j \neq i} x_j(k)}{|N_i(k)|}  \right| x(k), G(k) \right],  \ee
each node $j \neq i$ will be counted $\bar{m} \times \frac{m}{N-1}$ times for those neighborhood types such that $j \in N_i(k)$ and $|N_i(k)| = (m+1)$, so it holds that
\be \ba{lll} && \disp \sum_{G(k)} \left[ \left. \frac{\sum_{j\in N_i(k),j \neq i} x_j(k)}{|N_i(k)|}  \right| x(k), G(k) \right] \\ &=& \disp \sum_{j \neq i} \sum_{m = 0} ^ {N-1} \left[ \left. \frac{x_j(k)} {m+1} \bar{m} \frac{m}{N-1} \right| x(k) \right] . \ea \ee
Thus, we have
\be \ba{lll} & \E [ v_i(k) | x(k) ] = \displaystyle c_1 x_i(k)+
\sum_{j \neq i} c_2 x_j(k) ,\ea \ee
where
\be \ba{ll}  c_1  &:=  \displaystyle \frac{2^{(N-1)(N-2)/2} }{n}
\sum_{m=0} ^{N-1}
 {N-1 \choose m} \times \frac{1}{m+1},\\
c_2 
 & = \displaystyle \frac{2^{(N-1)(N-2)/2} }{n(N-1)}
 \sum_{m=1} ^{N-1}   {N-1 \choose m} \times
\frac{m}{m+1} . \ea \ee
This yields that, in view of (\ref{Exik}),
\be \E [ x_i(k+1) | x(k) ] = \frac{1}{\eta} \left[ c_1 x_i(k)
+ c_2 \sum_{j \neq i} x_j(k) \right], \ee
and hence
\be \ba{ll} &\E [ S(k+1) | x(k) ] \\
=& \displaystyle \sum_{i=1}^N
\E [ x_i(k+1) | x(k) ] \\
=& \displaystyle \frac{1}{\eta} \left[ c_1 S(k) + c_2(N-1) S(k) | x(k) \right] \\
=& \displaystyle \frac{1}{\eta} \frac{1}{2^{N-1}} 
         \sum_{m=0}^{N-1} {N-1 \choose m}  ( S(k)| x(k) )\\
         =&
         \displaystyle \frac{1}{\eta} \sum _{i=1}^N x_i(k). \ea \ee
Therefore, the expected state sum at the next time is
\be \ba{lll}  \E ( S(k+1) ) 
&=& \displaystyle \E [ \E
( S(k+1) | x(k) )]\\
&=& \displaystyle \frac{1}{\eta} \E    S(k) . \ea \ee
Since $\eta>1$, the above recursion converges to zero exponentially, and the decay
exponent is
\be -\frac{1}{k} \ln \frac{\E S(k)}{\E S(0)} =- \ln \frac{1}{\eta} = \ln \eta.\ee
\endproof \vsp

\textbf{Proof of Proposition \ref{prop:random}:} If $0<\eta \leq 1$, from Lemma
\ref{lemma:random1}, the system state sum will converge to the absorbing states with
probability 1.  If $\eta >1$, from Lemma \ref{lemma:random2}, the
system state will converge to zero exponentially with probability 1. \endproof  \vsp

\section{Numerical results}

\begin{figure}[h!]
\begin{center} \subfigure[]
{\epsfig{file=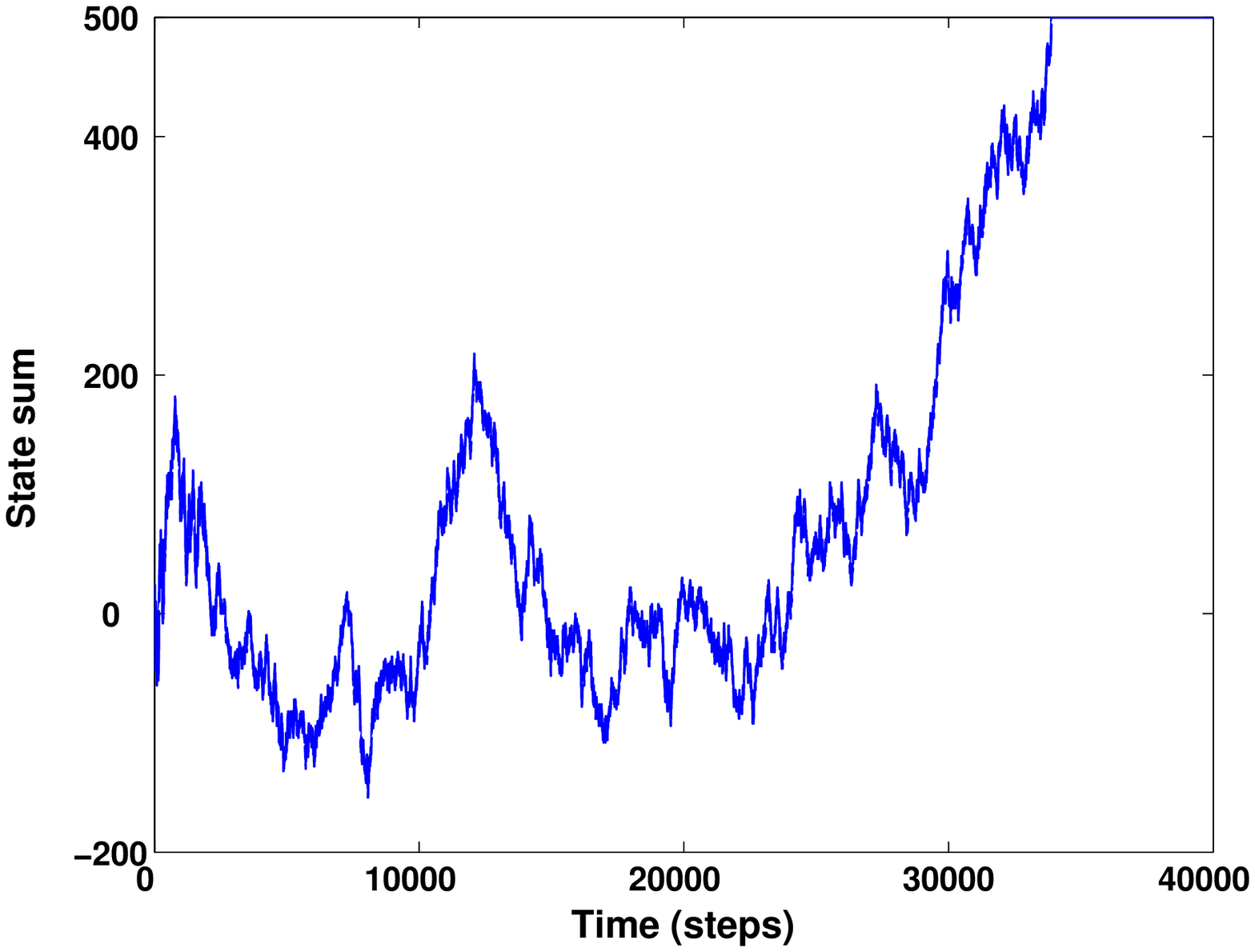,width=.48\hsize,height=0.36\hsize}
} \hspace{0pt} \subfigure[]
{\epsfig{file=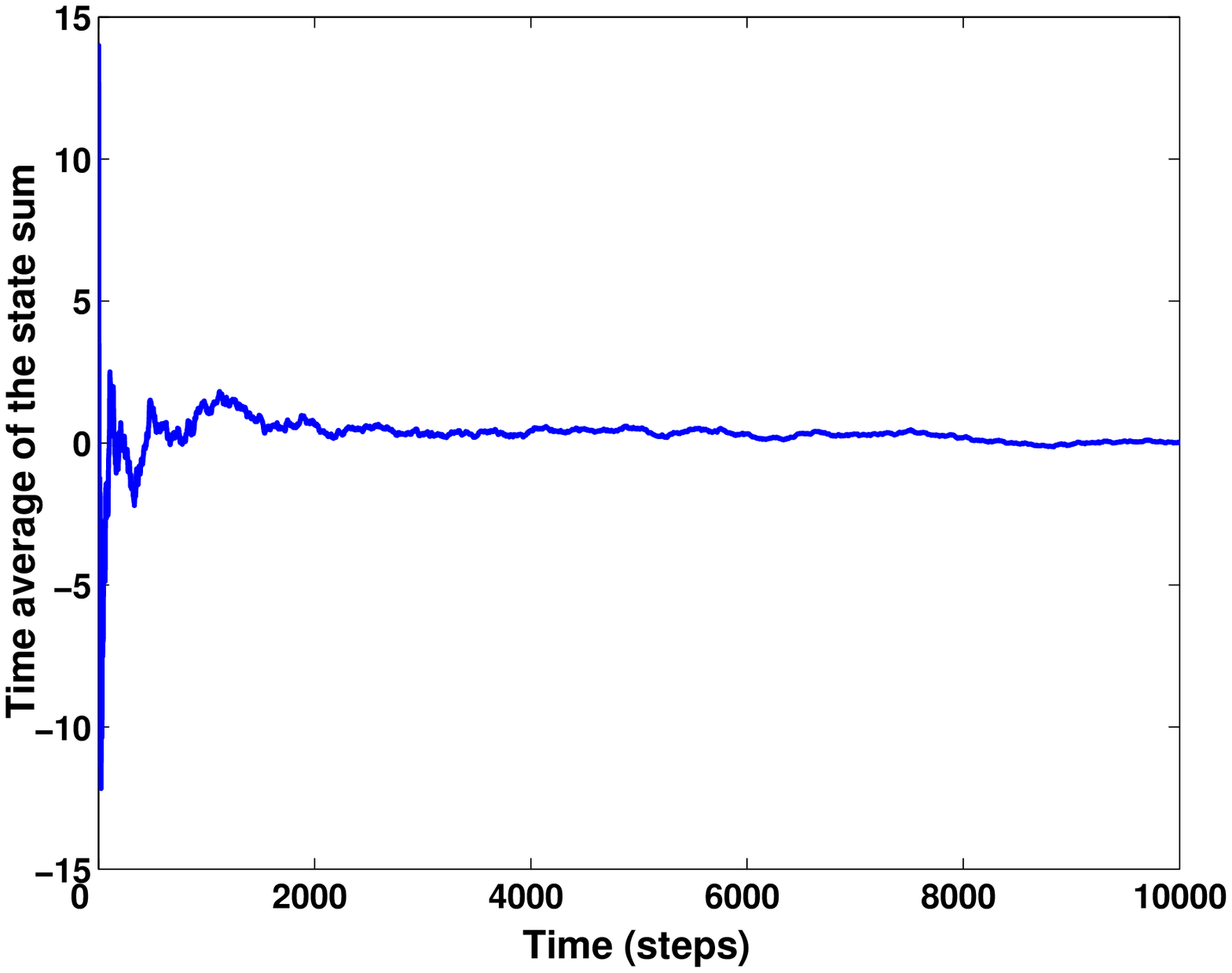,width=.47\hsize,height=0.35\hsize} }
\caption{Fixed graph simulation.  (a) Noise level is 0.75, and the
system converges to agreement of all $+1$s. (b) Noise level is 2,
and the system reaches disagreement in which about half of the
states are $+1$s and the other half are $-1$s.} \label{fig:fig1}
\end{center}
\end{figure}

Consider first a fixed one-dimensional 500-agent system. The agents are arranged along a
circle and each agent has two neighbors. The initial value of every agent is
arbitrarily assigned to be $+1$ or $-1$.  The simulation results demonstrate the phase
transition, see Figure \ref{fig:fig1} (a) and (b). In Figure \ref{fig:fig1} (a), the
vertical axis represents the state sum of the system, and the horizontal axis represents
the simulation steps. Figure \ref{fig:fig1} (a) demonstrates that, when the noise level is such
that $1/3<\eta \leq 1$, then all node values converge to agreement of all $+1$s (or all
$-1$s), that is, the state sum of the system is +500 (or -500). In Figure \ref{fig:fig1} (b),
the vertical axis represents the time average of the state sum, and the horizontal axis is
for the simulation steps. By ergodicity of the system, the time average should
converge to the ensemble average of the state sum. Figure \ref{fig:fig1} (b) shows that,
if the noise level is such that $\eta>1$, then the nodes reach
disagreement in which about half of the node values are $+1$s and the other half are
$-1$s.

\begin{figure}[h!]
\begin{center} \subfigure[]
{\epsfig{file=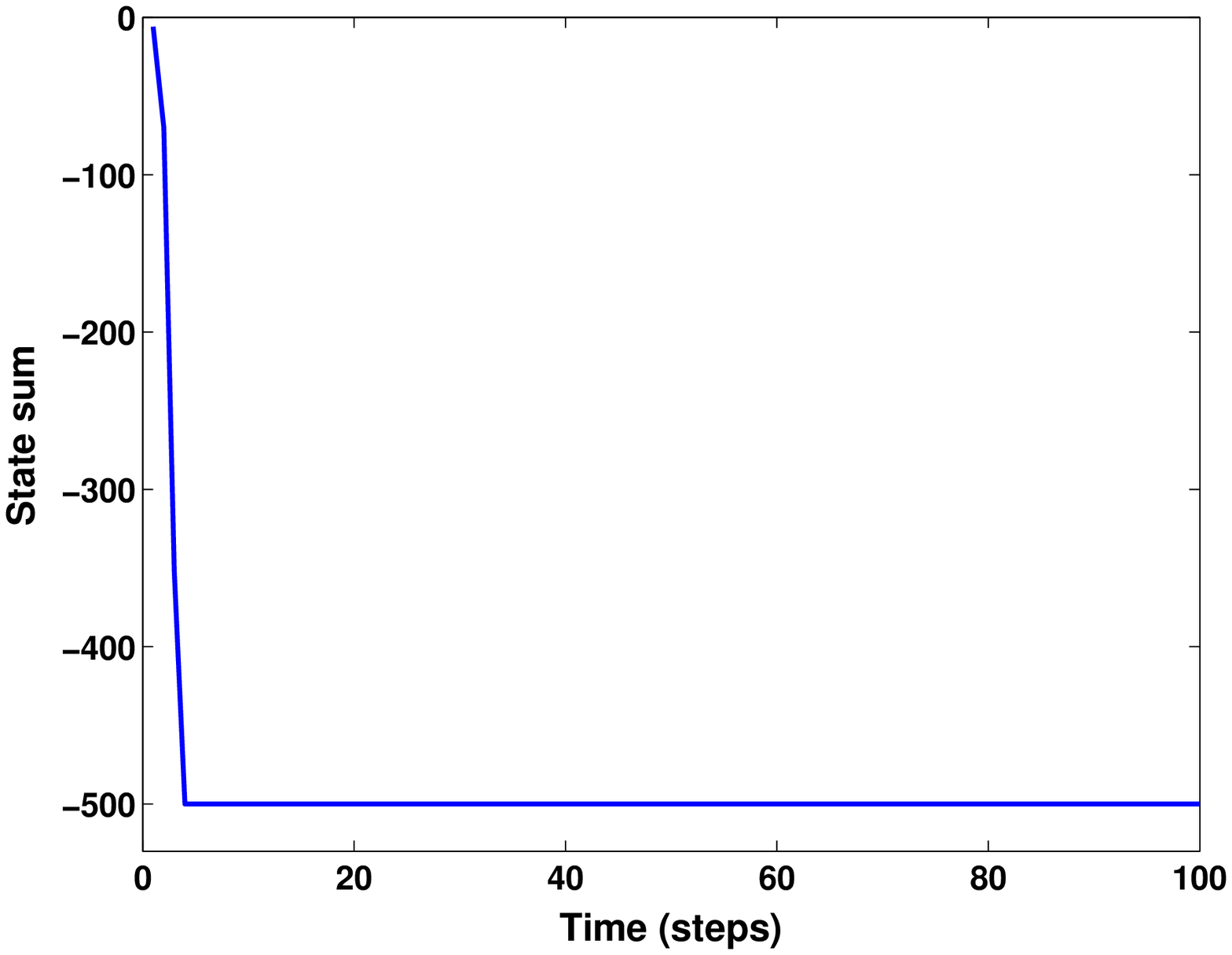,width=.48\hsize,height=0.36\hsize}
} \hspace{0pt}  \subfigure[]
{\epsfig{file=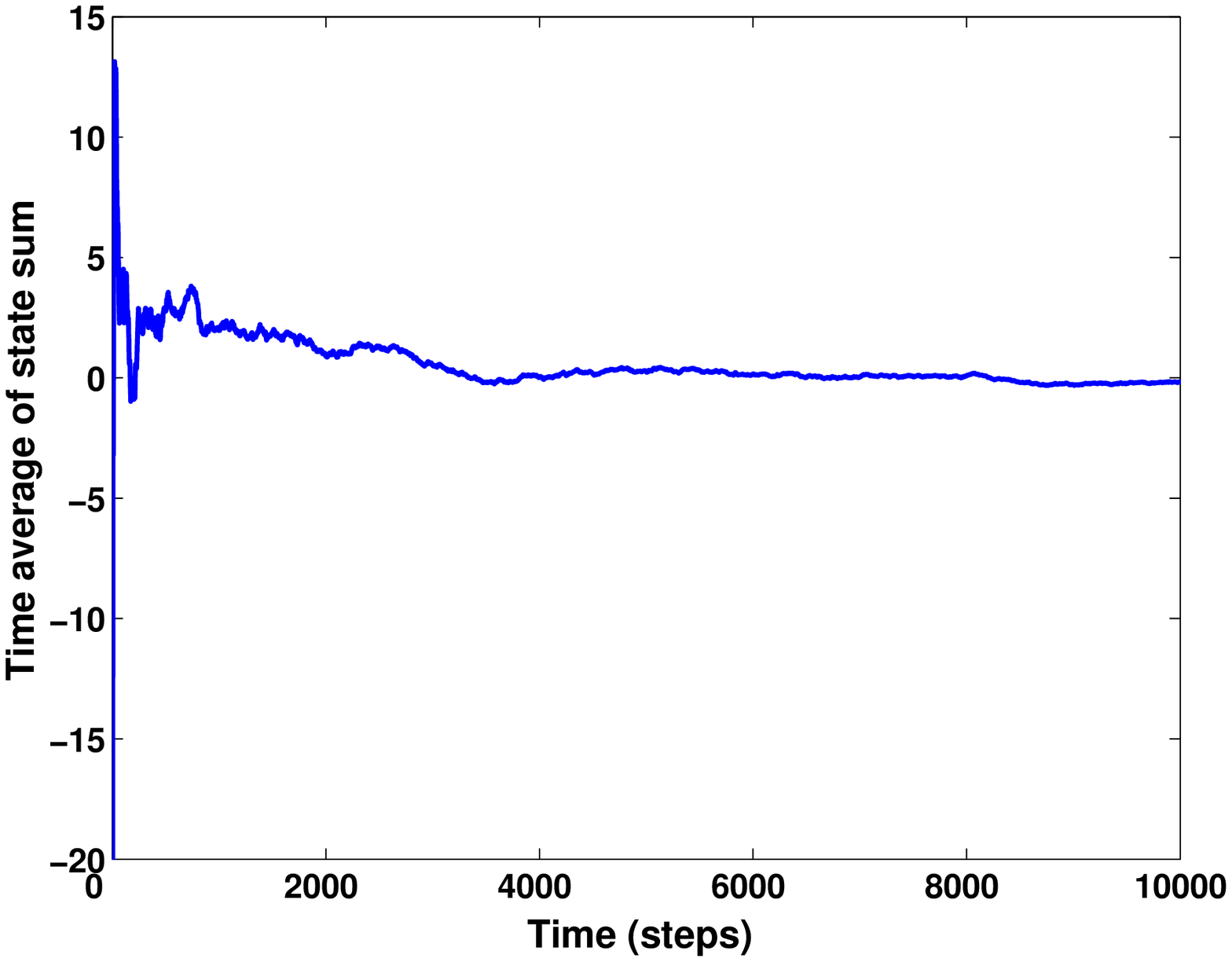,width=.47\hsize,height=0.35\hsize}
} \caption{Random graph process simulation. (a) Noise level is 0.005,
$p=0.1$, and the system converges to agreement of all $-1$s. (b)
Noise level is 2, $p=0.2$, and the
system reaches disagreement in which about half of the states
are $+1$s and the other half are $-1$s.} \label{fig:fig3}
\end{center}
\end{figure}

For the random graph process case, in our simulation we consider
\emph{binomial random graphs}. In a binomial random graph, each edge
has a probability $p$ to be formed at each time step and is
independent of all other edges and other times.  This means that
to generate such a binomial random graph, we only need to generate
at each step an adjacency matrix whose entries in the strictly upper
triangular part are independent and identically distributed.  The
initial value of every agent is randomly assigned to be $+1$ or
$-1$ according to an arbitrary distribution. The simulation
results are shown in Figure \ref{fig:fig3} (a) and (b), and are
similar to the fixed connected graph case, except that in the
random graph case, an arbitrarily small but positive noise level
can lead to agreement.

\begin{figure}[h!]
\begin{center}
{\scalebox{.6}
{\includegraphics{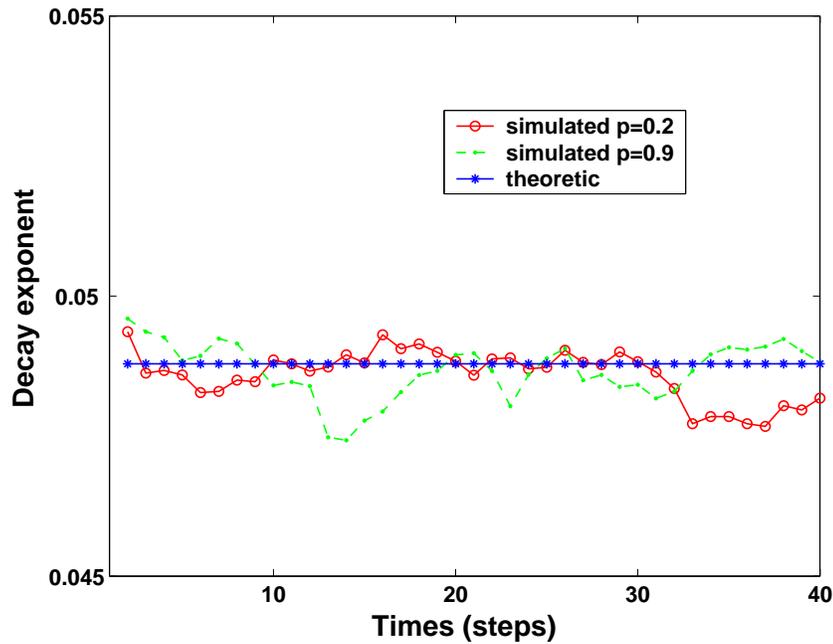}}} \caption{The simulated decay exponents
(averaged over 10,000 independent trials and $\eta=1.05$) and the theoretic decay
exponent.} \label{fig:exp}
\end{center}
\end{figure}

We can also compute the decay exponent of $\E S(k)$ from the numerical results. To obtain the probability means $\E S(k)$ numerically, we can run many independent trials of the random
process and take the average of the state sums across the trials.  See Figure \ref{fig:exp} for
the simulated decay exponents (with different edge probability $p$) and the theoretic decay exponent $\ln \eta$, which are
almost identical.  

Notice that $p$ does not play any role in the decay exponent. The role of $p$ is reflected in other quantities, such as the stationary distribution.  To see this, let us consider a two-node binomial graph, i.e. $N = 2$, and compute the stationary distribution as well as the decay exponent directly based on the transition probability matrix.  Let us first order the state values as $(+1,+1)$, $(+1,-1)$, $(-1,-1)$, and $(-1,+1)$.  Based on this ordering, the transition probability matrix is

\be P = \left( \ba{ c  c  c  c } 
		c 	& \disp \tabfrac{p}{4} + q b & 	a & \disp 	\tabfrac{p}{4} + q b \\ 
		b 	& \disp \tabfrac{p}{4} + q c & 	b & 	\disp \tabfrac{p}{4} + q a \\ 
		a 	& \disp \tabfrac{p}{4} + q b & 	c & 	\disp \tabfrac{p}{4} + q b \\ 
		b 	& \disp \tabfrac{p}{4} + q a & 	b & 	\disp \tabfrac{p}{4} + q c 
		\ea \right)
\ee
and the stationary distribution is
\be \ba{lllll}
\pi_{++}&=&\pi_{--}&=&\disp  \frac{ p + 4 q b}{ 2( p + 4 (1+q ) b )} \\
\pi_{+-}&=&\pi_{-+}&=&\disp  \frac{ 2 b}{ p + 4 (1+q ) b}, \ea
\ee%
%
where $q:= (1-p)$, and
%
\be 	a  :=  \disp \frac{( \eta - 1 ) ^2 }{4 \eta ^2}, 
		b  :=  \disp \frac{( \eta - 1 ) ( \eta + 1 ) }{4 \eta ^2} ,
		c  :=  \disp \frac{( \eta + 1 ) ^2 }{4 \eta ^2} .
\ee
Clearly, the edge probability $p$ influences the stationary distribution.  Now assume that the state is distributed according to distribution $p_0 := ( p_{++}, p_{+-}, p_{-+}, p_{--} )'$, which has the expected state sum as $2(p_{++} - p_{--} )$.  Then the state at the next time is distributed as $ P \times p_0$ and the expected state sum becomes $2(p_{++} - p_{--} )/\eta$.  Therefore, we have verified the dependence on $p$ in the stationary distribution and the independence on $p$ in the decay rate of the state sum.

\section{Conclusions and future work}

In this paper, we proposed simple dynamical systems models
exhibiting sharp phase transitions, and presented complete,
rigorous proofs of the phase transitions, with thresholds
found analytically.  Our analysis also provided a characterization
of how information (or noise) affects the collective behavior of
multi-agent systems, which gives an analytic explanation to the
intuition that, to reach consensus, high quality of communication
is needed.  These results hold for any dimension; in contrast,
phase transitions in the well known Ising models do not occur for
dimension one, and for dimension three or higher, Ising models are
NP complete and intractable.

In particular, we have shown that for a fixed connected graph, if
the noise level is greater than $(1 - 2/D)$ and less than $1$, all
the agents reach agreement, i.e. the state sum of the system
converges to $\pm N$, the only absorbing states of the system. For
noise level larger than $1$, the group of agents fail to reach any
agreement; instead they reach ``complete disagreement'' or disorder. Thus, a phase transition occurs at $\eta = 1$. For random
graph processes, the system reaches agreement even for noise level
smaller than $(1 - 2/D)$. This is because randomization is immune
to the artifacts (or local attractors) for smaller noise
which stops fixed graphs from reaching any agreement. However, the
tradeoff is that in random graphs, the nodes' neighbors may not be ``geographically
close", which might not hold true in some practical situations.

Our study was concentrated on the leaderless case. The leader case
is when there is a leader with a fixed value and it tries to
convince all other agents to follow its value. Simulation obtained
in this case suggested that a complete analysis is a bit involved
especially in the high noise regime, which is subject to further
research. Another direction could be to obtain a suitable Lyapunov
function for the models.  One advantage of doing so is that the
Lyapunov function based approach may be extended to rather general
nonlinear systems, as suggested by \cite{meg:lyap,lucMoreau05}.
The Lyapunov function is preferably a quadratic one, leading to mean-squre
stability, which is stronger than the
mean stability obtained in this paper.  The
applications of our approach and results are also subject to
future research, including the extension of our approach to more
realistic models; note that our models in this paper are simple
and not quite realistic, though the simplicity helped us to
completely characterize the phase transition.  We will also
explore the connections of our model to relevant models, e.g.
the Ising models, Hopfield networks, cellular automata, other
random graphs, etc.  Finally, we remark that the approach and results developed in this paper may be found useful to study more general dynamical systems under communication constraints, such as cooperation with limited communication, complex systems in the presence of noise, etc.  The study of such problems would help establish insights on how information (or limited information) interacts with system dynamics to generate various types of interesting system behavior.

\appendix

We prove that $\pi(x) = \pi(-x)$ for any $x$ in four steps. 

Step 1: Establish a one-to-one mapping between the $2J$ possible values (see (\ref{ss_size})) that the state of the system can take and integers $\pm 1, \pm 2, \cdots, \pm J$, such that if state $x$ is mapped to $+j$, then state $-x$ is mapped to $-j$.  Now aggregate the states as follows.  Let $\bar{j}:=
(j,-j)$ for any $j = 1, \cdots, J$.  Then we induce from the Markov process $\{x(k)\} _{k=0} ^\infty$ another Markov process $\{ \bar{x}(k) \} _{k=0} ^\infty$, where the latter is defined on the induced state space consisting of all $\bar{j}$s.  Note that it is straightforward to verify that
$\{ \bar{x}_k \}  _{k=0} ^\infty$ forms a Markov
process on the induced state space, and this Markov process is ergodic.

Step 2: Denote the transition probability matrix for process $\{ \bar{x}(k)\}  _{k=0} ^\infty $ as
$\bar{p}$, and the corresponding stationary distribution vector as
$\bar{\pi}
:=(\bar{\pi}(\bar{1}), \bar{\pi}(\bar{2}), ...,
\bar{\pi}(\bar{J}))'
$. Then it holds that $\bar{\pi} = \bar{p}
\bar{\pi}$. By ergodicity, $\bar{\pi}$ is non-zero and unique
(i.e., the matrix $(I-\bar{p})$ must be rank deficient).

Step 3: For the Markov process $\{x(k)\} _{k=0} ^\infty$, denote the stationary distribution
vector as $\pi:=(\pi_1',\pi_2')'$, where $\pi_1
=(\pi(+1),\pi(+2),..., \pi(+J))'$ and $\pi_2
=(\pi(-1),\pi(-2),..., \pi(-J))'$.  It can be verified that, by
the symmetry that the state transition $i \rightarrow j$ has the same probability as the state transition $(-i) \rightarrow (-j)$, the transition probability
matrix has the following particular form:
\be p:=\left( \ba{cc} A & B \\ B & A \ea \right). \ee

Step 4: By ergodicity of $\{x(k)\} _{k=0} ^\infty$, it holds that 
\be \pi = p \pi \label{erg:short} \ee
or equivalently,
\be \ba{lll} \pi_1 &=& A \pi_1 + B \pi_2 \\
\pi_2 &=& B \pi_1 + A \pi_2 . \ea \label{eq:erg}\ee
However, it can be easily seen that $\bar{p} = A+B$.  Notice that
\be \pi_1 = \pi_2 = \bar{\pi} \label{sol}\ee
solves (\ref{eq:erg}), i.e., $\pi_0 := ( \bar{\pi}', \bar{\pi}' )'$ solves (\ref{erg:short}) and is non-zero.  By ergodicity, the
non-zero solution is unique, and hence $\pi_0$ must be the
solution to (\ref{erg:short}), which follows that $\pi(j) = \pi(-j)$ for any $j$ or $\pi(x) = \pi(-x)$ for any $x$.

\vsp \vsp
 \small\hfill \markright{\textsf{References}}  \vsp

\bibliographystyle{unsrt}

\end{document}